\documentclass[reqno]{amsart}

\usepackage{style}

\author{Mireille Boutin}
\address{School of Electrical and Computer Engineering, Purdue
  University, West Lafayette, IN, USA}
\email{mboutin@purdue.edu}
\author{Gregor Kemper} \address{Technische Universit\"at M\"unchen,
  Zentrum Mathematik - M11, Boltzmannstr. 3, 85748 Garching, Germany}
\email{kemper@ma.tum.de}

\title{A Drone Can Hear the Shape of a Room}

\date{\today}


\begin{document}

\begin{abstract}
  We show that one can reconstruct the shape of a room with planar
  walls from the first-order echoes received by four non-planar
  microphones placed on a drone with generic position and orientation.
  Both the cases where the source is located in the room and on the
  drone are considered. If the microphone positions are picked at
  random, then with probability one, the location of any wall is
  correctly reconstructed as long as it is heard by four
  microphones. Our algorithm uses a simple echo sorting criterion to
  recover the wall assignments for the echoes. We prove that, if the
  position and orientation of the drone on which the microphones are
  mounted do not lie on a certain set of dimension at most 5 in the
  6-dimensional space of all drone positions and orientations, then
  the wall assignment obtained through our echo sorting criterion must
  be the right one and thus the reconstruction obtained through our
  algorithm is correct. Our proof uses methods from computational
  commutative algebra.
\end{abstract}

\maketitle

\section*{Introduction} \label{sIntro}%

Assume we have a room, by which we understand an arrangement of planar
walls, which may include ceilings, floors, and sloping walls.  An
omnidirectional loudspeaker and some omnidirectional microphones are
in the room.  The loudspeaker, modeled as a point source, emits a
short duration pressure wave (a sound impulse) at a frequency high
enough so that the ray acoustics approximation is valid.  The
microphones receive several delayed responses corresponding to the
sound bouncing back from each wall. These are the first-order
echoes. These echoes subsequently bounce back from each wall again,
creating second-order echoes, and so on.  We are interested in the
problem of reconstructing the shape of the room from the first-order
echoes.  Specifically, we use the time delay of each first-order
echoes, in other words the propagation time, which provides us with a
set of distances from each microphone to mirror images of the source
reflected across each wall. Since we do not know which echo
corresponds to which wall, the distances are unlabeled.  In fact,
depending on the microphone configuration and room geometry, a
microphone may not receive any echo from a given wall.  The problem is
to figure out under which circumstances, and how, one can find out the
correct distance-wall assignments and reconstruct the wall positions.

The distances to the mirror images of the source are obtained from the
time of arrivals of the impulses at each microphone. If the sound
impulse is known, and if the microphones and loudspeaker share a
common clock, these times of arrival can be computed by finding the
relevant peaks of the cross-correlations between the original impulse
and the signals received by the microphones.
Since the sound impulse gets more and more blurred as it bounces from
wall to wall, it is possible to distinguish the peaks corresponding to
the first order echoes from those of the higher order ones, assuming
that the signal to noise ratio is large enough.  See
\mycite{antonacci2012inference} for a discussion on how this can be
accomplished and the practical limitation of these assumptions.  We
shall assume that there are no missed peaks and no spurious peaks. We
shall also assume that the peaks corresponding to first-order echoes
can be distinguished from those of higher-order echoes.

This work focuses on theoretically determining when the problem is
well-posed (i.e., when the distances contain enough information to
uniquely reconstruct the room) in the minimal case of four microphones
and one source. Aside from \mycite{DPWLV1}, most other authors have
set aside such theoretical questions and worked on developing
numerical solution methods.  Seeking robustness, the problem is often
set up experimentally so to create redundancies and/or simplifying
assumptions.  For example, a direct numerical solution method in which
one searches for the wall positions on a discrete grid was proposed by
\mycite{CroccoICASSP2018}.  In that work, the number of walls is
assumed to be known, the room is assumed to be a convex polyhedron,
and several sources and microphones are used in order to increase the
robustness of the reconstruction.  Different problem definitions are
considered by other authors. For example, a single microphone acquires
the data of a source moving in a circle around it in the work of
\mycite{AntonacciICASSP2010}. The case where the microphone positions
are unknown [\citenumber{Thrun2005}] is also of interest.  A slightly
different formulation where one or more microphones are moving along
an unknown path in a room with unknown geometry containing one or more
(potentially moving) speakers is considered, a problem called echo
SLAM (simultaneous localization and mapping), for example
[\citenumber{hu2011simultaneous},\citenumber{krekovic2016echoslam},\citenumber{krekovic2019shapes}].
Other times, the known geometry of the room is used to localize an
indoor object carrying some receivers through the multiple reflections
of a signal bouncing off the walls (multipath propagation), for
example
[\citenumber{MeissnerThesis2014},\citenumber{leitinger2015evaluation},\citenumber{witrisal2016high},\citenumber{shahmansoori2017position},\citenumber{mendrzik2018harnessing}].

The core of our paper is \cref{1aDetect}, which describes a procedure
to detect walls from first-order echoes acquired by four microphones
whose positions are known. When we think of a wall, we distinguish
between a wall and the plane in which is is contained. In particular,
walls are usually finite, as is illustrated in~\cref{1fBad}. By {\em
  detecting} a wall we mean that four non-collinear points on the wall
are determined. Clearly this uniquely determines the plane containing
the wall, but also provides some information about the actual location
of the wall within that plane.  The key part of the algorithm is a
simple echo sorting criterion (Relation~\cref{1eqRelation}) that is
used to solve the wall assignment problem. The criterion is a
vanishing Cayley-Menger determinant involving the pairwise distances
between the five-point configuration formed by the four microphones
and the mirror image of the source through one wall. As the criterion
is always satisfied by the distances corresponding to a correct wall
assignment, the algorithm detects all walls that are heard by the four
microphones. However, it can sometimes detect walls that are not there
({\em ghost walls}).

Our main results are \cref{2tFixed,3tMain}, which specify conditions
under which \cref{1aDetect} is guaranteed not to detect any wall that
is not there (no ghost wall). \Cref{2tFixed} assumes that the
loudspeaker is at a fixed location inside the room. \Cref{3tMain}
assumes that the loudspeaker is carried by the drone, along with the
microphones. The conditions for both theorems are very general; they
are satisfied with probability one if the orientation and position of
the drone are picked at random following a non-degenerate probability
density function. Specifically, the set of exceptional drone positions
and orientations lie inside a subvariety of dimension at most five
within the six-dimensional space of possible drone placements.

\mycite{DPWLV1} have considered the case where the microphones and the
loudspeaker are placed in a fixed location inside the room.  For
reasons of simplicity, the room was assumed to be a convex polyhedron,
and the two reconstruction methods presented used five
microphones. The authors' Theorem~1, which applies to four or more
microphones, guarantees the correctness of the reconstruction for all
but a set of measure zero of microphone arrangements. Notice that
in~[\citenumber{DPWLV1}], the microphone arrangements are chosen
randomly from a 12-dimensional configuration space, whereas for our
result a six-dimensional space suffices. Also, placing the microphones
on a drone rather than independently in the room opens up new
application scenarios.  We show that our \cref{2tFixed} implies
Theorem~1 of~[\citenumber{DPWLV1}] with \cref{2cDPWLV}.

The key to proving \cref{2tFixed} is \cref{2Claim2}, which states that
one can rotate and translate the drone so it lies in a position for
which it is guaranteed that our echo sorting criterion will lead to
correct wall assignments. The proof of this claim is accomplished with
the help of the symbolic computation software
MAGMA~[\citenumber{magma}]. \Cref{3tMain} is proved in a similar
fashion. \\

%
{\bf Acknowledgments.} We thank Ivan Dokmani{\'{c}} for fruitful
conversations. Our thanks also go to the anonymous referees for their
careful reading of the manuscript and for many comments that helped us
improve the paper.

\section{Preliminaries}
\label{sPreliminaries}

\newcommand{\ve}[1]{\mathbf{#1}}%

\begin{figure}[htbp]
  \newcommand{\vertex}[1][black]{\node[#1,circle, draw, outer sep=2pt,
    inner sep=0pt, minimum size=3pt,fill]}
  \centering
  \begin{tikzpicture}[scale=1.0]
    \vertex (L) at (0,0) [label={[black]left:$\ve L$}] {};
    \vertex (s1) at (0,4)  [color=red,label={right:$\ve s_1$}] {};
    \vertex (s2) at (6,0)  [color=blue,label={right:$\ve s_2$}] {};
    \vertex (ss) at (-6,0)  [color=white,label={[white]left:$ss$}] {};
    \vertex (m1) at (-1,1)  [label={[black]135:$\ve m_1$}] {};
    \vertex (m2) at (1,1)  [label={[black]45:$\ve m_2$}] {};
    \vertex (m3) at (-1,-1)  [label={[black]-135:$\ve m_3$}] {};
    \vertex (m4) at (1,-1)  [label={[black]-45:$\ve m_4$}] {};
    
    \draw[red,ultra thick] (-3,2) -- (3,2);
    \draw (-2.7,1.9) node[below]{$W_1$};
    \fill[pattern=north west lines,pattern color=red] (-3,2) --
    (-3,1.1*2) -- (3.2,1.1*2) -- (3,2) -- (-3,2);
    \draw[color=red,very thick,dotted] (L) -- (2/3,2);
    \draw[color=red,very thick,dotted,>=stealth,->] (2/3,2) -- (m2);
    \draw[color=red,dashed] (2/3,2) -- (s1);
    \draw[color=red,very thick,dotted] (L) -- (-2/3,2);
    \draw[color=red,very thick,dotted,>=stealth,->] (-2/3,2) -- (m1);
    \draw[color=red,dashed] (-2/3,2) -- (s1);
    \draw[color=red,very thick,dotted] (L) -- (2/5,2);
    \draw[color=red,very thick,dotted,>=stealth,->] (2/5,2) -- (m4);
    \draw[color=red,dashed] (2/5,2) -- (s1);
    \draw[color=red,very thick,dotted] (L) -- (-2/5,2);
    \draw[color=red,very thick,dotted,>=stealth,->] (-2/5,2) -- (m3);
    \draw[color=red,dashed] (-2/5,2) -- (s1);

    \draw[blue,ultra thick] (3,2) -- (3,-2);
    \draw (3.3,-1.8) node[right]{$W_2$};
    \fill[pattern=north west lines,pattern color=blue] (3,2) --
    (3.2,1.1*2) -- (3.2,-2) -- (3,-2) -- (3,2);
    \draw[color=blue,very thick,dotted] (L) -- (3,3/5);
    \draw[color=blue,very thick,dotted,>=stealth,->] (3,3/5) -- (m2);
    \draw[color=blue,dashed] (3,3/5) -- (s2);
    \draw[color=blue,very thick,dotted] (L) -- (3,-3/5);
    \draw[color=blue,very thick,dotted,>=stealth,->] (3,-3/5) -- (m4);
    \draw[color=blue,dashed] (3,-3/5) -- (s2);
    \draw[color=blue,very thick,dotted] (L) -- (3,3/7);
    \draw[color=blue,very thick,dotted,>=stealth,->] (3,3/7) -- (m1);
    \draw[color=blue,dashed] (3,3/7) -- (s2);
    \draw[color=blue,very thick,dotted] (L) -- (3,-3/7);
    \draw[color=blue,very thick,dotted,>=stealth,->] (3,-3/7) -- (m3);
    \draw[color=blue,dashed] (3,-3/7) -- (s2);
  \end{tikzpicture} \\ \mbox{} \\
  \caption{Each microphone hears the echoes from two walls. Virtually,
  the sound comes from the mirror points~$\ve s_1$ and~$\ve s_2$.}
  \label{1fMirror}
\end{figure}
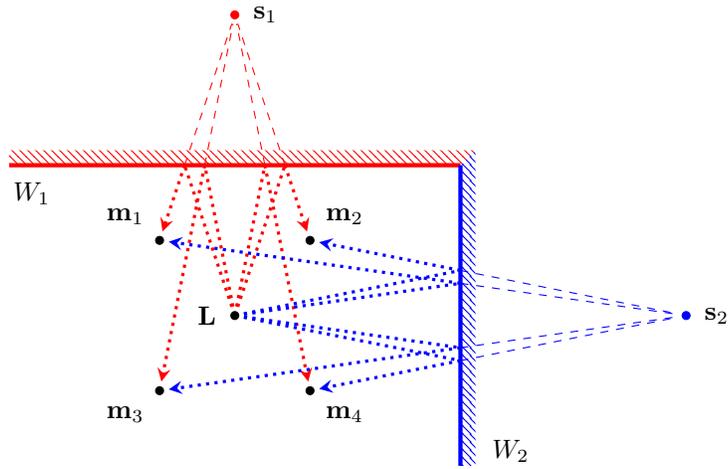

For what we have to say in this section, it is irrelevant whether the
loudspeaker or the microphones are placed at fixed positions or
mounted on a drone. Following the image method of
\mycite{allen1979image}, let $\ve s \in \RR^3$ be the point obtained
by reflecting the loudspeaker position with respect to one of the
walls (or, more precisely, the plane containing the wall). We call
such an~$\ve s$ a \df{mirror point}. A sound emitted from the
loudspeaker and reflected at the wall and then heard by the
microphones corresponds to a sound emitted from the mirror
point~$\ve s$ and traveling directly to the microphones. So by
measuring the time elapsed between sound emission and echo detection,
the distance $\lVert\ve s - \ve m_i\rVert$ from~$\ve{s}$ to the
microphone at position~$\ve m_i$ can be determined. This is
illustrated in \cref{1fMirror}.

The following proposition states some facts about four microphones
hearing the echo from one wall. Part~\ref{1pDetectD} gives a relation
between the squared distances between~$\ve s$ and the microphones. The
relation is just a restatement the well-known fact that the
Cayley-Menger determinant of five points in $\RR^3$ vanishes (see
\mycite{Cayley:1841}).

\begin{prop} \label{1pDetect}%
  In the above situation, assume that four microphones at positions
  $\ve m_i \in \RR^3$ ($i = 1 \upto 4$) hear the sound reflected at a
  wall (the same wall for all microphones).
  \begin{enumerate}[label=(\alph*)]
  \item \label{1pDetectA}%
    With
    \begin{equation} \label{1eqM}%
      M := \left(
        \begin{array}{cccc}
          \\
          \ve m_1 & \ve m_2 & \ve m_3 & \ve m_4 \\ \\
          \hline %
          1 & 1 & 1 & 1
        \end{array}\right) \in \RR^{4 \times 4},
    \end{equation}
    the microphones are coplanar if and only if $\det(M) = 0$.
  \item \label{1pDetectB}%
    Assume from now on that the microphones are non-coplanar and write
    $\tilde{M} \in \RR^{3 \times 4}$ for the upper $3 \times 4$-part
    of $(M^{-1})^T$, the transpose inverse of $M$. Then~$\ve s$ can be
    computed from the squared distances $d_i := \lVert\ve s - \ve
    m_i\lVert^2$ by
    \begin{equation} \label{1eqComputeS}%
      \ve s = \frac{1}{2} \tilde{M} \cdot
      \begin{pmatrix}
        \lVert\ve m_1\rVert^2 - d_1 \\
        \vdots \\
        \lVert\ve m_4\rVert^2 - d_4
      \end{pmatrix}.
    \end{equation}
  \item \label{1pDetectC}%
    Let $\ve L \in \RR^3$ be the position of the loudspeaker. Then the
    wall at which the sound was reflected lies on the plane with
    normal vector $\ve s - \ve L$ and passing through the point
    $\frac{1}{2}(\ve s + \ve L)$. Four non-collinear points on the wall
    can be found by intersecting the line between~$\ve s$ and
    $\ve m_i$ ($1 \le i \le 4$) with this plane. These points are
    given by
      \[
      (1 - \tau_i) \ve s + \tau_i \ve m_i \quad \text{with} \quad
      \tau_i = \frac{\lVert \ve s - \ve L\rVert^2}{2 \langle \ve s -
        \ve L,\ve s - \ve m_i\rangle},
      \]
      where $\langle\cdot,\cdot\rangle$ denotes the standard scalar
      product.
    \item \label{1pDetectD}%
      With $D_{i,j} := \lVert\ve m_i - \ve m_j\rVert^2$ and
      $u_1 \upto u_4 \in \RR$ any numbers, set
      \begin{equation} \label{1eqF}%
        D :=
        \begin{pmatrix}
          0 & u_1 & \cdots & u_4 & 1 \\
          u_1 & D_{1,1} & \cdots & D_{1,4} & 1 \\
          \vdots & \vdots & & \vdots & \vdots \\
          u_4 & D_{4,1} & \cdots & D_{4,4} & 1 \\
          1 & 1 & \cdots & 1 & 0
        \end{pmatrix} \in \RR^{6 \times 6} \quad \text{and} \quad
        f_M(u_1 \upto u_4) := \det(D).
      \end{equation}
      Then the~$d_i$ from~\ref{1pDetectB} satisfy the relation
      \begin{equation} \label{1eqRelation}%
        f_M(d_1 \upto d_4) = 0.
      \end{equation}
  \end{enumerate}
\end{prop}

Before giving the proof, we consider an example of the relation~$f_M$
from~\ref{1pDetectD}.

\begin{ex} \label{1exStandard}%
  Consider the configuration of microphones given by
  \[
  M = \left(
  \begin{array}{cccc}
    0 & 0 & 1 & 0 \\
    0 & 0 & 0 & 1 \\
    0 & 1 & 0 & 0 \\
    \hline
    1 & 1 & 1 & 1
  \end{array}\right)
  \]
  (so the microphones are at the origin and the standard basis
  vectors). For this we have
  \[
  f_M(u_1 \upto u_4) = 4 (u_2 - u_1 - 1)^2 + 4 (u_3 - u_1 - 1)^2 + 4
  (u_4 - u_1 - 1)^2 - 16 u_1.
  \]
  %
  %
  %
  Since the coefficients of $f_M$ only depend on the relative
  distances between the microphones, they do not change when the
  microphones are moved while maintaining their relative position.
\end{ex}

\begin{proof}[Proof of \cref{1pDetect}]
  The equivalence in~\ref{1pDetectA} can be obtained by subtracting
  the first column of $M$ from the other columns. Part~\ref{1pDetectB}
  follows from
  \begin{equation} \label{eqMsequalDistances}%
    M^T \cdot \left(
      \begin{array}{c}
        \\ 2 \ve s \\ \\ \hline -\lVert\ve s\rVert^2
      \end{array}
    \right) =
    \begin{pmatrix}
      2 \langle\ve m_1,\ve s\rangle - \lVert\ve s\rVert^2 \\
      \vdots \\
      2 \langle\ve m_4,\ve s\rangle - \lVert\ve s\rVert^2
    \end{pmatrix} =
    \begin{pmatrix}
      \lVert\ve m_1\rVert^2 - \lVert\ve m_1 - \ve s\rVert^2 \\
      \vdots \\
      \lVert\ve m_4\rVert^2 - \lVert\ve m_4 - \ve s\rVert^2
    \end{pmatrix} = \begin{pmatrix}
      \lVert\ve m_1\rVert^2 - d_1 \\
      \vdots \\
      \lVert\ve m_4\rVert^2 - d_4
    \end{pmatrix}.
  \end{equation}
  The first statement from~\ref{1pDetectC} follows from the
  definition of~$\ve s$. It follows from ray acoustics that the
  intersection of the line $\overline{\ve s\ve m_i}$ with the plane
  lies on the wall, and from the non-coplanarity of the~$\ve m_i$ that
  the intersection points are non-collinear. The equation for~$\tau_i$
  follows from verifying that the
  $\ve x_i := (1 - \tau_i) \ve s + \tau_i \ve m_i$ satisfy the
  equation
  $\langle\ve x_i,\ve s - \ve L\rangle = \frac{1}{2}\langle \ve s +
  \ve L,\ve s - \ve L\rangle$
  defining the plane. (One may also check that $0 < \tau_i < 1$ holds
  if and only if~$\ve m_i$ and~$\ve L$ lie on the same side of the
  plane, which must be true if the echo is heard by the $i$th
  microphone.)

  For the proof of~\ref{1pDetectD} we refer to \mycite{Cayley:1841}.
\end{proof}

\begin{rem} \label{1rL}%
  Assume that the positions $\ve m_i$ of the microphones are known,
  but the position $\ve L$ of the loudspeaker is not. Then $\ve L$ can
  be determined, using \cref{1pDetect}\ref{1pDetectB}, as follows: The
  sound traveling directly from the loudspeaker to the microphones can
  always be distinguished from the echoes since it is the first to
  arrive at the microphones. Therefore the distances between $\ve L$
  and the $\ve m_i$ can be determined. So applying
  \cref{1pDetect}\ref{1pDetectB} to these distances yields $\ve L$.
\end{rem}

\section{Remarks about Relation~\cref{1eqRelation}}

Relation~\cref{1eqRelation} is really the relation, given in
\mycite[Proposition~2.2(b)]{Boutin.Kemper}, satisfied by the pairwise
distances between five points in $3$-space. For example, using the
forth microphone as an anchor point (point $n$ in
$\Delta_{ij}:=d_{in}+d_{jn}-d_{ij}$ where $d_{ij}$ is the squared
distance from point i to point j) among the four microphones
$\ve m_1, \ve m_2, \ve m_3, \ve m_4$, and setting $\ve m_0:=\ve s$,
the matrix $\ve \Delta=(\Delta_{ij})$ becomes
\begin{eqnarray*}
  \ve  \Delta&=&
                 \left(   \| \ve m_i-\ve m_4 \|^2 +\|\ve m_j-\ve
                 m_4\|^2-\|\ve m_i-\ve m_j\|^2 \right)_{i,j=0,1,2,3}, \\
             &=&  \left(  D_{i,4} +D_{j,4} - D_{i,j} \right)_{i,j=0,1,2,3}.
\end{eqnarray*}
Its determinant is zero because it can be factored as
\begin{equation} \label{eq:factorization}%
  \ve \Delta= 2 \left( (\ve m_i-\ve m_4 )^T (\ve m_j-\ve m_4)
  \right)_{i,j=0,1,2,3},
\end{equation}
and therefore has rank at most three. Replacing $D_{i,0}$ in the
matrix $\ve \Delta$ by an indeterminate $u_i$, for $i=1,2,3,4$, we
obtain the following relation, which is equivalent to
Relation~\cref{1eqRelation} up to a sign:
 \begin{equation}
 \label{2eqRelation}
\begin{small}
 \det 
 \left( 
  \begin{array}{cccc}
2 u_4 & u_4+D_{1,4}-u_1 & u_4+D_{2,4}-u_2 & u_4+D_{3,4}-u_3\\
 u_4+D_{1,4}-u_1 &  2 D_{1,4}  & D_{1,4}+D_{2,4}-D_{1,2} &D_{1,4}+D_{3,4}-D_{1,3} \\
 u_4+D_{2,4}-u_2 &   D_{2,4}+D_{1,4}-D_{1,2} & 2 D_{2,4} & D_{2,4}+D_{3,4}-D_{2,3} \\
 u_4+D_{3,4} -u_3 &  D_{1,4}+D_{3,4}-D_{1,3} &  D_{2,4}+D_{3,4}-D_{2,3} & 2D_{3,4}
 \end{array}
 \right)
 =0.
\end{small}
 \end{equation}
Neither Relation~\cref{1eqRelation} nor Relation~\cref{2eqRelation}
is sufficient to guarantee the existence of a point configuration with
the corresponding distances.  However if the matrix $\Delta$ is
positive semi-definite, then its eigendecomposition
$\Delta= Q^T\Lambda Q$, where $\Lambda$ is a diagonal matrix with at
most three non-zero diagonal elements, yields a solution
\begin{equation}
\sqrt{2}
  \left(  \begin{array}{cccc}
\ve m_0'-\ve m_4' & \ve m_1'-\ve m_4' & \ve m_2'-\ve m_4' & \ve m_3'-\ve m_4'  \\ 
0 & 0 & 0 & 0
\end{array} \right)
=
\sqrt{ \Lambda}  Q.
\label{eq:svd_solution}  
\end{equation}
If there are three non-zero eigenvalues of $\ve \Delta$ and they are
distinct, then it can be shown that this is the unique solution, up to
an orthogonal transform, for the factorization of $\Delta$
as~\cref{eq:factorization}. In other words, there exists a rotation
mapping the correct solution $\ve m_i- \ve m_4$ to the points
$\ve m_i'- \ve m_4'$ obtained by eigenvalue decomposition
from~\cref{eq:svd_solution}.  Assuming that the positions of the
microphones $\ve m_1, \ve m_2, \ve m_3, \ve m_4$ are known, the
rotation can be computed as the product
\[
 \left(  \begin{array}{cccc}
\ve m_1'-\ve m_4' & \ve m_2'-\ve m_4' & \ve m_3'-\ve m_4'  
\end{array} \right)  \left(  \begin{array}{cccc}
\ve m_1-\ve m_4 & \ve m_2-\ve m_4 & \ve m_3-\ve m_4  
\end{array} \right)^{-1}. 
\]

Observe that this provides an alternative to
Relation~\cref{1eqRelation}, namely that the matrix $\Delta$ (or,
equivalently, the matrix $D$ evaluated at $u_i=d_i$) has at most 3
non-zero eigenvalues. When the distance measurements are inaccurate,
the fourth eigenvalue could become non-zero. If the error is small, an
approximate reconstruction could be obtained by setting the smallest
eigenvalue to zero. However in the vicinity of a bad drone position,
the conditioning of this reconstruction method could be very bad, as
it would be impossible to distinguish between small perturbations of
the different reconstructions possible for the bad drone position.

Note that both Relation~\cref{1eqRelation} and
Relation~\cref{2eqRelation} are different from the echo sorting
criteria used in the reconstruction methods
of~[\citenumber{DPWLV1}]. The two criteria used in the reconstruction
both use at least 5 microphones.  One criterion, derived from
Equation~\cref{eqMsequalDistances}, is that the range of $M$ should
include the vector $(\|\ve m_1\|^2-d_1, \ldots, \|\ve
m_5\|^2-d_5)^T$. The other criterion, which is the one used in the
proposed reconstruction method of [\citenumber{DPWLV1}], uses the {\em
  Euclidean distance matrix} (EDM). Specifically, it uses the 6-by-6
matrix of pairwise distances $\left( D_{i,j}\right)$ between the five
(or more) microphones and the mirror point $\ve s$. The criterion is
that the rank of this EDM matrix is at most five.

In the practical algorithm proposed in [\citenumber{DPWLV1}] to handle
noisy data, the classical (metric) multidimensional scaling technique
[\citenumber{abdi2007metric}] is used to reconstruct the points from
the approximate distance measurements. This is done by first {\em
  centering} the Euclidean distance matrix
$D_{EDM}=\left( D_{i,j}\right)$ as
\[
E=-B D_{EDM} B 
\]
where $B=I_{N+1}-\frac{1}{ (N+1)} \left(1,1,\ldots,1 \right) \left( 1,1,\ldots,1\right)^T $ 
and $N+1$ is the total number of points used ($N$ microphones plus one wall mirror point $\ve s=\ve m_0$). The matrix $E$ is called the {\em inner product matrix}; indeed, if we map the center of mass of the points $\ve m_i$ to the origin before taking their inner product, we get 
\begin{multline*}
2  \left( \ve m_i-\sum_{k=0}^N \frac{\ve m_k}{N+1} \right)^T \left( \ve m_j-\sum_{k=0}^N \frac{\ve m_k}{N+1} \right) = \\
  -D_{i,j}+\sum_{k=0}^N \frac{D_{k,i}}{N+1}+\sum_{m=0}^N\frac{D_{m,j}}{N+1}- \sum_{k,m=0}^N \frac{D_{k,m}}{(N+1)^2}
= E_{i,j}.
\end{multline*}
Note that this factorization is valid for any number of
microphones. Thus, the rank of the centered matrix $E$ is at most
three, for any number of microphones, and a solution (up to an
orthogonal transform) can be obtained by eigendecomposition of $E$.

%

\section{The wall detection algorithm} 
\label{3sDetection}

We now give an algorithm that attempts to detect walls from the echoes
heard by four microphones. Since each microphone hears echoes from
multiple walls, it is necessary to match those echoes that come from
the same wall. A natural way do this is to use the
relation~\cref{1eqRelation}. Notice that by ``detecting'' a wall we
mean that four non-collinear points on the wall are computed; no
further information about the actual expanse of the wall within the
plane containing it can be obtained from the echoes of a single sound
emission, unless additional hypotheses are made on the walls,
e.g. that they are the facets of a convex polyhedron (see
\mycite{DPWLV1}).

\begin{algorithm}
  \caption{Detect walls from first-order echoes} \label{1aDetect}
  \mbox{}%
  \begin{description}
  \item[Input] The delay times of the first-order echoes recorded by
    four microphones, and the distances $D_{i,j}$ between the
    microphones.
  \end{description}
  \begin{algorithmic}[1]
    \STATE \label{1aDetect1} For $i = 1 \upto 4$, collect the delay
    times of the first-order echoes recorded by the $i$th microphone
    in the set $\mathcal T_i$.%
    \STATE \label{1aDetect2} Set
    $\mathcal D_i := \{c^2 (t - t_0)^2 \mid t \in \mathcal T_i\}$
    ($i = 1 \upto 4$), where~$c$ is the speed of sound and~$t_0$ is
    the time of sound emission.%
    \FOR{$(d_1,d_2,d_3,d_4) \in \mathcal D_1 \times \mathcal D_2 \times
      \mathcal D_3 \times \mathcal D_4$} \label{1aDetect3}%
    \STATE \label{1aDetect4} With~$f_M$ defined by~\cref{1eqF},
    evaluate $f_M(d_1 \upto d_4)$.%
    \IF{$f_M(d_1 \upto d_4) = 0$}%
    \STATE \label{1aDetect5} Use~\cref{1eqComputeS} to compute the
    mirror point~$\ve s$ from $(d_1 \upto d_4)$.%
    \STATE \label{1aDetect6} Use \cref{1pDetect}\ref{1pDetectC} to
    compute four non-collinear points on the wall with mirror
    point~$\ve s$ and, if desired, a normal vector.%
    \STATE \label{1aDetect7} {\bf Output} the data of this wall.%
    \ENDIF%
    \ENDFOR%
  \end{algorithmic}
\end{algorithm}

If for
$(d_1 \upto d_4) \in \mathcal D_1 \times \cdots \times \mathcal D_4$,
the $d_i$ come from echoes from the {\em same} wall, then the relation
$f_M(d_1 \upto d_4) = 0$ holds and therefore the wall will be
detected. So the algorithm is guaranteed to detect every wall from
which a first-order echo is heard by all microphones. It is possible,
however, that the algorithm detects walls that are not really there
(``ghost walls''; see \cref{1exGhost}). The main purpose of this
paper is to show that these mistakes are rare.

Note that the search for matches in step~\ref{1aDetect3} can be
accelerated by using the triangle inequalities. For example, one could
start with the shortest distances, say $d_1$, for the first
microphone. Then the distance to the second microphone $d_2$ would
need to satisfy the inequality
$\sqrt{d_2}\leq \sqrt{d_1}+\|\ve m_1-\ve m_2 \| $, and so on.

The following example shows that it can happen that the algorithm
detects ghost walls.

\begin{ex} \label{1exGhost}%
  \cref{1fBad}%
  \newcommand{\vertex}[1][black]{\node[#1,circle, draw, outer sep=3pt,
    inner sep=0pt, minimum size=3pt,fill]}
  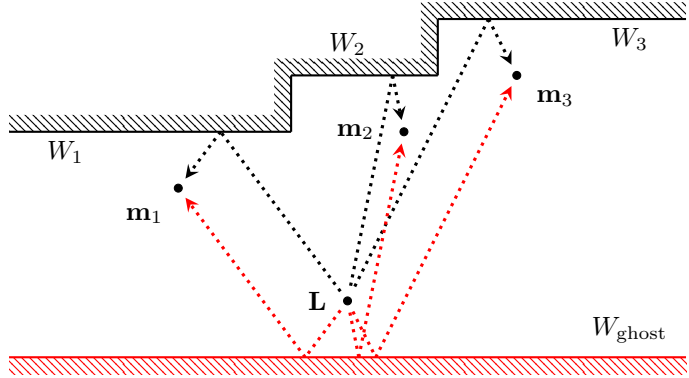
\begin{figure}[htbp]
    \centering
    \begin{tikzpicture}[scale=0.75]
      \fill[pattern=north west lines,pattern color=red] (1,1) --
      (13,1) -- (13,0.7) -- (1,0.7) -- (1,1);
      \draw[red,thick] (1,1) -- (13,1);
      \draw (12,1.1) node[above]{$W_{\operatorname{ghost}}$};
      \vertex (L) at (7,2)  [label=left:$\ve L$] {};
      \vertex (m1) at (4,4) [label=-120:$\ve m_1$] {};
      \vertex (m2) at (8,5) [label=left:$\ve m_2\ $] {};
      \vertex (m3) at (10,6) [label=-30:$\ve m_3$] {};
      \fill[pattern=north west lines] (1,5) -- (6,5) -- (6,6) --
      (8.6,6) -- (8.6,7) -- (13,7) -- (13,7.3) -- (8.3,7.3) --
      (8.3,6.3) -- (5.7,6.3) -- (5.7,5.3) -- (1,5.3) -- (1,5);
      \draw[thick] (1,5) -- (6,5);
      \draw (2,5) node[below]{$W_1$};
      \draw[thick] (6,5) -- (6,6);
      \draw[thick] (6,6) -- (8.6,6);
      \draw (7,6.2) node[above]{$W_2$};
      \draw[thick] (8.6,6) -- (8.6,7);
      \draw[thick] (8.6,7) -- (13,7);
      \draw (12,7) node[below]{$W_3$};
      \draw[very thick,dotted,>=stealth,->] (3/4*4+7/4,5) -- (m1);
      \draw[very thick,dotted] (L) -- (3/4*4+7/4,5);
      \draw[very thick,dotted,>=stealth,->] (4/5*8+1/5*7,6) -- (m2);
      \draw[very thick,dotted] (4/5*8+1/5*7,6) -- (L);
      \draw[very thick,dotted,>=stealth,->] (5/6*10+7/6,7) -- (m3);
      \draw[very thick,dotted] (5/6*10+7/6,7) -- (L);
      \draw[very thick,dotted,red,>=stealth,->] (3/4*7+4/4,1) -- (m1);
      \draw[very thick,dotted,red] (L) -- (3/4*7+4/4,1);
      \draw[very thick,dotted,red,>=stealth,->] (4/5*7+8/5,1) -- (m2);
      \draw[very thick,dotted,red] (L) -- (4/5*7+8/5,1);
      \draw[very thick,dotted,red,>=stealth,->] (5/6*7+10/6,1) -- (m3);
      \draw[very thick,dotted,red] (L) -- (5/6*7+10/6,1);
    \end{tikzpicture}
    \caption{The microphones $\ve m_i$ think they are hearing echoes
      from the wall $W_{\operatorname{ghost}}$. The dotted lines
      stand for sound rays, $\ve L$ for the loudspeaker.}
    \label{1fBad}
  \end{figure}
  shows three microphones in a plane at positions~$\ve m_1$,
  $\ve m_2$, $\ve m_3$ that hear echoes from three walls~$W_i$, but
  the time elapsed between sound emission and echo detection is the
  same as if they were hearing echoes from one single wall
  $W_{\operatorname{ghost}}$, which does not exist. This arises
  because (1) the walls $W_i$ are all parallel to each other, (2)
  each~$\ve m_i$ can hear the echo from $W_i$, and (3) the distances
  between $\ve m_i$ and $W_i$ are the same for all~$i$. It is easy to
  add a fourth microphone outside of the plane, together with a wall
  possibly also outside of the plane, such that~(1)--(3) extend to the
  fourth microphone and wall. (We find it harder to include that in
  our two-dimensional sketch in \cref{1fBad}.) Since the echoes
  heard by the microphones {\em could} have come from the single wall
  $W_{\operatorname{ghost}}$, the relation~\cref{1eqRelation} is
  satisfied, and so \cref{1aDetect} will falsely detect
  $W_{\operatorname{ghost}}$ as a wall.
  %
\end{ex}

One can argue that such bad examples correspond to exceptional wall
configurations.  Indeed, it is easy to show that, for wall
configurations picked at random following a non-degenerate probability
distribution, there is a probability zero of this happening.
Specifically, there is a probability zero of the wall configuration
yielding distances to the microphones that satisfy
$f_M(\| \ve s_1-\ve m_1 \|^2, \| \ve s_2-\ve m_2 \|^2, \| \ve s_3-\ve
m_3 \|^2, \| \ve s_4-\ve m_4 \|^2)=0$ with the $\ve s_i$'s not all
equal.  For example, suppose that the first point $\ve s_1$ is the
correct mirror point $\ve s$ and that another point say
$\ve s_4\neq \ve s$, but $\ve s_1=\ve s_2=\ve s_3=\ve s$.  Then one
can freely change the value of the distance $\| \ve s_4- \ve m_4\|^2$
by moving the wall corresponding to $\ve s_4$.  The distance
$\| \ve s_4- \ve m_4\|^2$ is then changed freely without affecting the
other distances
$\| \ve s_1-\ve m_1 \|^2, \| \ve s_2-\ve m_2 \|^2, \| \ve s_3-\ve m_3
\|^2$ because these correspond to another wall, and thus the zero set
of $f_M$ can be avoided. Similarly, if the first two points
$\ve s_1=\ve s_2= \ve s$ correspond to the correct wall and the last
two points $\ve s_3=\ve s_4\neq \ve s$ correspond to another wall,
then the values of the distances
$\| \ve s_3-\ve m_3 \|^2, \| \ve s_4-\ve m_4 \|^2$ can be changed
freely by moving to wall corresponding to $\ve s_3$ and $\ve s_4$ in
${\mathbb R}^3$. This will not affect the first two distances
$\| \ve s_1-\ve m_1 \|^2, \| \ve s_2-\ve m_2 \|^2$ and thus the zero
set of $f_M$ can be avoided. Note that our argument holds because
moving the walls does not affect the relationship
$f_M(u_1,u_2,u_3,u_4)=0$ itself, since the coefficients of its
polynomial are monomials in the distances between the
microphones. Therefore, we can conclude that all but a set of measure
zero of choices of mirror points $\ve s_1, \ve s_2, \ve s_3, \ve s_4$
that are not all equal satisfy
$f_M(\| \ve s_1-\ve m_1 \|^2, \| \ve s_2-\ve m_2 \|^2, \| \ve s_3-\ve
m_3 \|^2, \| \ve s_4-\ve m_4 \|^2)\neq 0$.  For fixed microphone
positions $\ve m_1, \ve m_2, \ve m_3, \ve m_4$, one could thus make
sure that the wall configuration is not an exceptional one by checking
that the mirror points satisfy the
inequalities: 
\[
f_M(\|\ve s_1-\ve m_1 \|^2, \| \ve s_2-\ve m_2 \|^2, \| \ve s_3-\ve m_3 \|^2, \| \ve s_4-\ve m_4 \|^2)\neq 0, \text{ for all } \{\ve s_1, \ve s_2, \ve s_3, \ve s_4\} \text{ not all equal}.
\]
However, moving the walls is not practical, and it is conceivable that the exceptional wall configurations might include all cases that have, say, parallel walls. 
Therefore we seek to show instead that, given any wall configuration, most microphone positions satisfy the above inequalities.
The argument required for this is more difficult. 

\section{A drone in a room with a fixed loudspeaker}
\label{2sFixed}

As in the previous section, assume that we are given a room with a
loudspeaker and four microphones in it. In this section we will assume
that the loudspeaker is at a fixed position, but the microphones are
mounted on a drone (or, mathematically speaking, that their relative
positions are fixed). 
Apart from the
fact that this is a realistic scenario for applications, this has the
computational advantage that the coefficients of the relation
$f_M(d_1 \upto d_4) = 0$ that is exploited in the wall detection
algorithm remain the same once and for all, since by
\cref{1pDetect}\ref{1pDetectD} the coefficients only depend on the
mutual distances between the microphones. For applying
\cref{1aDetect}, we need to assume that all the microphone locations
are known, which may be problematic when the microphones are on a
drone. Two answers can be given to this objection:
\begin{enumerate}[label=\arabic*.]
\item In some scenarios, it can really be assumed that the drone knows
  its position and orientation.  For example, the drone can be equipped 
  with an inertial measurement unit or its position can be calibrated using
  an external positioning system based on triangulation.
\item If the drone is not aware of its own position and orientation,
  it can perform the computation with respect to ``its own'' coordinate
  system. More precisely, this means that the coordinates of the
  microphone positions $\ve m_i$ are assigned once and for all,
  according to where the microphones are located on the
  drone. \cref{1aDetect} will then detect the walls with respect to
  the coordinate system, traveling with the drone, in which the
  microphone coordinates were assigned. Step~\ref{1aDetect6} of the
  algorithm requires the position $\ve L$ of the loudspeaker. This,
  too, can be measured within the drone's coordinate system by using
  \cref{1rL}.
\end{enumerate}

We say that the microphones (or the drone) are in a \df{good position}
if \cref{1aDetect} detects no walls that are not really there (no
ghost wall.)  Recall that the algorithm is guaranteed to detect every
wall from which a first-order echo is heard by all microphones.
Therefore if the drone is moving from a bad position to a good
position, all the ghost walls previously reconstructed will disappear
and all the walls that are there and can still be heard by the
microphones will remain. If additional walls can be heard from the new
position, then these will be added to the set of reconstructed walls
as well.
The goal of this section is to prove the following result, which
implies that  generic drone positions are good positions.

\begin{theorem} \label{2tFixed}%
  Consider a given room, by which we understand an arrangement of
  walls, which may include ceilings, floors, and sloping walls. Assume
  there is a loudspeaker at a given position in the room. Also
  consider a drone that carries four non-coplanar microphones at fixed
  locations on its body. Place the drone in the room at a random
  position, which means that not only the location of the drone's
  center of gravity is chosen at random, but also its pitch, yaw, and
  roll. Then with probability~$1$ the drone is in a good
  position. More precisely, within the configuration space
  $\RR^3 \times \SO(3)$ of possible drone positions, the bad ones lie
  in a subvariety of dimension~$\le 5$.
\end{theorem}

As a consequence we obtain the main result of \mycite{DPWLV1}.

\begin{cor} \label{2cDPWLV}%
  We make the same assumptions as in \cref{2tFixed}, except that the
  microphones are not mounted on a drone, but are placed independently
  at random locations. Then with probability~$1$ they are in a good
  position. More precisely, within the configuration space $\RR^{12}$
  of possible microphone positions, the bad ones lie in a subvariety of
  dimension~$\le 11$.
\end{cor}

Since it is not immediately clear how the corollary follows from the
theorem, both will be proved together. Before giving the actual proof,
we present a rough roadmap, aiming to convey the geometry that lies
behind the proof.

\begin{enumerate}[label=\arabic*.]
\item We define a subset $\mathcal U \subseteq \RR^3 \times \SO(3)$ of
  so-called ``very good drone positions.'' The set $\mathcal U$
  depends on the given arrangement of walls and on the given
  configuration of microphones, and turns out to be Zariski open.
\item We reduce to the case of four walls (see \cref{2Claim2}
  below).
\item After suitable choices of coordinates, we can model an
  arrangement of four walls together with a configuration of four
  microphones by a single point
  $(b_1 \upto b_6,c_1 \upto c_5) \in \RR^{11}$. To express the
  dependencies, let us write
  ${\mathcal U}(b_1 \upto b_6,c_1 \upto c_5)$ instead of $\mathcal
  U$. Notice that a point in $\RR^{11}$ may encode a wall arrangement
  where all or some walls are equal.
\item We consider the set
  \[
    \mathcal V := \bigl\{(b_1 \upto b_6,c_1 \upto c_5) \in \RR^{11} \mid
    {\mathcal U}(b_1 \upto b_6,c_1 \upto c_5) = \emptyset\bigr\}
    \subseteq \RR^{11}
  \]
  of wall arrangements and microphone configurations for which there
  is no very good drone position. On the other hand, consider the set
  $\mathcal V' \subseteq \RR^{11}$ of all
  $(b_1 \upto b_6,c_1 \upto c_5)$ such that $(c_1 \upto c_5)$ defines
  a coplanar microphone configuration, or the walls given by
  $(b_1 \upto b_6)$ are all equal. Notice that we are now allowing the
  wall arrangement and microphone configuration to vary.
\item Now it is enough to show that
  \begin{equation} \label{2eqGeometry}%
    \mathcal V \subseteq \mathcal V'.
  \end{equation}
  Indeed, if a wall arrangement and a microphone configuration satisfy
  the hypothesis of \cref{2tFixed}, then the corresponding point does
  not lie in $\mathcal V'$ and therefore, by~\cref{2eqGeometry}, also
  not in $\mathcal V$. This means
  ${\mathcal U}(b_1 \upto b_6,c_1 \upto c_5) \ne \emptyset$, so by the
  Zariski openness, the set of bad drone positions lies in a proper
  subvariety of the configuration space.
\item The inclusion~\cref{2eqGeometry} translates to the
  ideal-theoretic statement~\cref{2eqJ} below. In that statement, the
  ideal $J$ defines the set $\mathcal V$. (The sets $\mathcal V$ and
  $\mathcal V'$ are never mentioned in the proof. We only need them
  here to explain the underlying geometry.)
\item Finally,~\cref{2eqJ} is verified by Gr\"obner basis
  computations, using a computer algebra system. The computation is
  demanding because the ideal $J$ has many generators. As detailed in
  the proof, some tricks contribute to the feasibility of the
  computation.
\end{enumerate}

\newcommand{\ini}{{\operatorname{ini}}}%
\newcommand{\Mini}{M_\ini}%
\newcommand{\Aff}{\operatorname{Aff}}%
\newcommand{\ASO}{\operatorname{ASO}}%
\newcommand{\re}{\operatorname{ref}}%
\newcommand{\h}{{\operatorname{hom}}}
\begin{proof}[Proof of \cref{2tFixed,2cDPWLV}]
  In \cref{2tFixed}, the random placement of the drone means that from
  initial positions of the microphones, given by the~$\ve m_i^\ini$ in
  the matrix
  \begin{equation} \label{2eqMini}%
    \Mini = \left(
      \begin{array}{cccc}
        \\
        \ve m_1^\ini & \ve m_2^\ini & \ve m_3^\ini & \ve m_4^\ini \\
        \\
        \hline 1 & 1 & 1 & 1
      \end{array}\right) \in \RR^{4 \times 4},
  \end{equation}
  the actual positions~$\ve m_i$ are given by
  \begin{equation} \label{2eqM}%
    M := \left(
      \begin{array}{cccc}
        \\
        \ve m_1 & \ve m_2 & \ve m_3 & \ve m_4 \\ \\
        \hline %
        1 & 1 & 1 & 1
      \end{array}\right) = A \cdot \Mini \quad \text{with} \quad 
    A = \left(
      \begin{array}{cccc}
        a_{1,1} & a_{1,2} & a_{1,3} & a_{1,4} \\
        a_{2,1} & a_{2,2} & a_{2,3} & a_{2,4} \\
        a_{3,1} & a_{3,2} & a_{3,3} & a_{3,4} \\
        \hline 0 & 0 & 0 & 1
      \end{array}
    \right),
  \end{equation}
  where the upper left $3 \times 3$-part of $A$ lies in $\SO(3)$. In
  \cref{2cDPWLV}, by contrast, the matrix $A$ can be chosen freely
  without the $\SO(3)$-condition. Let us write $\Aff(3)$ for the
  $12$-dimensional space of all matrices $A$ as in~\cref{2eqM}, and
  $\ASO(3)$ for the subset where the upper left $3 \times 3$-part
  comes from $\SO(3)$. The configuration space $\RR^3 \times \SO(3)$
  can be identified with $\ASO(3)$. We call a matrix $A \in \Aff(3)$
  \df{good} if the microphone positions given by~\cref{2eqM} are
  good.

  Let $\mathcal W$ be the set of walls from our room. In this proof we
  identify the walls with the planes containing them. The mirror
  points are given by $\ve s = \re_W(\ve L)$, the reflection of the
  loudspeaker position at a wall $W \in \mathcal W$.  We call
  $A \in \Aff(3)$ \df{very good} if the following holds: For four
  walls $W_1 \upto W_4 \in \mathcal W$ the relation
  \[
  f_M\bigl(\lVert\re_{W_1}(\ve L) - \ve m_1\rVert^2 \upto
  \lVert\re_{W_4}(\ve L) - \ve m_4\rVert^2\bigr) = 0
  \]
  is satisfied only if $W_1 = W_2 = W_3 = W_4$. (Recall that the above
  expression depends on $A$ by~\cref{2eqM}.)
  \begin{claim} \label{2Claim1}
    If $A$ is very good, then it is good.
  \end{claim}
  Indeed, with the notation of \cref{1aDetect}, let
  $(d_1 \upto d_4) \in \mathcal D_1 \times \cdots \times \mathcal
  D_4$. For each~$i$ there exists a wall $W_i \in \mathcal W$ such
  that $d_i = \lVert\re_{W_i}(\ve L) - \ve m_i\rVert^2$. If
  $f_M(d_1 \upto d_4) = 0$, then $W_1 = W_2 = W_3 = W_4 =: W$ by
  hypothesis, so by \cref{1pDetect}, the wall that is rendered in
  step~\ref{1aDetect7} of the algorithm is $W$. This shows that the
  algorithm detects no ghost walls.
  \begin{claim} \label{2Claim2}%
    Let $\ve s_1,\ve s_2,\ve s_3,\ve s_4 \in \RR^3$ be four points
    that are not all equal. Then there exists $A \in \ASO(3)$ such
    that
    \[
    f_M\bigl(\lVert\ve s_1 - \ve m_1\rVert^2 \upto \lVert\ve s_4 - \ve
    m_4\rVert^2\bigr) \ne 0.
    \]
    (Recall that the above expression depends on $A$ by~\cref{2eqM}.)
  \end{claim}
  Before proving the claim, we show that it implies
  \cref{2tFixed,2cDPWLV}. By \cref{1pDetect}\ref{1pDetectC}, a wall
  $W \in \mathcal W$ is uniquely determined by its mirror point
  $\re_W(\ve L)$. Therefore the claim implies that for
  $W_1 \upto W_4 \in \mathcal W$ which are not all equal, the set
  \[
  \mathcal U_{W_1 \upto W_4} := \bigl\{A \in \ASO(3) \mid
  f_M\bigl(\lVert\re_{W_1}(\ve L) - \ve m_1\rVert^2 \upto
  \lVert\re_{W_4}(\ve L) - \ve m_4\rVert^2\bigr) \ne 0\bigr\}
  \]
  is non-empty. Since
  $f_M\bigl(\lVert\re_{W_1}(\ve L) - \ve m_1\rVert^2 \upto
  \lVert\re_{W_4}(\ve L) - \ve m_4\rVert^2\bigr)$ depends polynomially
  on the coefficients of $A$ and since $\ASO(3)$ is an irreducible
  variety, this implies that the complement of
  $\mathcal U_{W_1 \upto W_4}$ in $\ASO(3)$ has dimension strictly
  less than $\dim(\ASO(3)) = 6$. It follows that also the finite
  intersection
  \[
    \mathcal U := \bigcap_{\substack{W_1 \upto W_4 \in \mathcal W\
        \text{such that} \\
        \text{not all} \ W_i\ \text{are equal}}} \mathcal U_{W_1 \upto
      W_4}
  \]
  has a complement of dimension $\le 5$. By definition, all
  $A \in \mathcal U$ are very good and therefore, by
  \cref{2Claim1}, also good. So indeed every bad $A$ lies in the
  complement of $\mathcal U$, which is a subvariety of $\ASO(3)$ of
  dimension $\le 5$, and \cref{2tFixed} follows. To conclude
  \cref{2cDPWLV} one just needs to replace $\ASO(3)$ by $\Aff(3)$ in
  this argument.

  It remains to prove \cref{2Claim2}, for which we can forget
  about the arrangement of walls. So we are only given four vectors
  $\ve s_1,\ve s_2,\ve s_3,\ve s_4 \in \RR^3$, not all equal, and the
  matrix~$\Mini$ (see~\cref{2eqMini}).

  For the feasibility of the computation in the final part of the
  proof it is necessary to choose a convenient Cartesian coordinate
  system. First, $\ve s_1$ can be chosen as the origin of the
  coordinate system, so the matrix
  $S = \begin{pmatrix} \ve s_1 \ \ve s_2 \ \ve s_3 \ \ve
    s_4 \end{pmatrix} \in \RR^{3 \times 4}$
  becomes
  $S = \begin{pmatrix}\ve 0 \ \ve s_2 \ \ve s_3 \ \ve
    s_4 \end{pmatrix}$.
  Moreover, using QR-decomposition, we can write
  \[
  \begin{pmatrix}
    \ve s_2 \ \ve s_3 \ \ve s_4
  \end{pmatrix} = Q \cdot
  \begin{pmatrix}
    b_1 & b_2 & b_3 \\
    0 & b_4 & b_5 \\
    0 & 0 & b_6
  \end{pmatrix}
  \]
  with $Q \in \SO(3)$, $b_i \in \RR$. So if we use the columns of $Q$
  (instead of the standard basis vectors) as basis of $\RR^3$, $S$
  becomes
  \begin{equation} \label{2eqS}%
    S := \begin{pmatrix} \ve s_1 \ \ve s_2 \ \ve s_3 \ \ve s_4
    \end{pmatrix} =
    \begin{pmatrix}
      0 &  b_1 & b_2 & b_3 \\
      0 & 0 & b_4 & b_5 \\
      0 & 0 & 0 & b_6
    \end{pmatrix}.
  \end{equation}
  With this, the hypothesis that the $\ve s_i$ are not all equal
  becomes $S \ne 0$.

  We also need to simplify the matrix $\Mini$, given
  by~\cref{2eqMini}. Since we wish to prove \cref{2Claim2},
  which states that there exists $A \in \ASO(3)$ such that
  $A \cdot \Mini$ satisfies a certain condition, we may modify $\Mini$
  by multiplying it with suitable matrices from $\ASO(3)$ on the
  left. Writing
  $\Mini = \left(\begin{array}{cccc}\ve c_1 & \ve c_2 & \ve c_3 & \ve
      c_4 \\ \hline 1 & 1 & 1 & 1 \end{array}\right)$, we have
  $\left(\begin{array}{c|c} I_3 & -\ve c_1 \\ \hline 0 &
      1 \end{array}\right) \cdot \Mini = \left(\begin{array}{c|c}0 &
      C' \\ \hline 1 & 1 \end{array}\right)$ with
  $C' \in \RR^{3 \times 3}$. By QR-decomposition, we can write
  $C' = Q' R$ with $Q' \in \SO(3)$ and $R$ upper triangular, so
  \[
  \left(
    \begin{array}{c|c}
      (Q')^{-1} & 0 \\
      \hline 0 & 1
    \end{array}\right) \left(
    \begin{array}{c|c}
      0 & C' \\
      \hline 1 & 1
    \end{array}\right) = \left(
    \begin{array}{c|c}
      0 & R \\
      \hline 1 & 1
    \end{array} \right) =  \left(
    \begin{array}{cccc}
      0 & c_0 & c_1 & c_2 \\
      0 & 0 & c_3 & c_4 \\
      0 & 0 & 0 & c_5 \\
      \hline 1 & 1 & 1 & 1
    \end{array}\right)
  \]
  with $c_0 \upto c_5 \in \RR$. By \cref{1pDetect}\ref{1pDetectA},
  the hypothesis that the microphones are non-coplanar translates to
  $\det(\Mini) \ne 0$, so $c_0 c_3 c_5 \ne 0$. We can now rescale our
  coordinate system by a factor of~$c_0$ (which corresponds to the
  choice of a unit of length). In summary, we may assume
  \begin{equation} \label{2eqNewMini}%
    \Mini = \left(
      \begin{array}{cccc}
        0 & 1 & c_1 & c_2 \\
        0 & 0 & c_3 & c_4 \\
        0 & 0 & 0 & c_5 \\
        \hline 1 & 1 & 1 & 1
      \end{array}\right)
  \end{equation}
  with $c_i \in \RR$, $c_3 c_5 \ne 0$.

  Turning to the proof of \cref{2Claim2}, we observe that
  $f_M\bigl(\lVert\ve s_1 - \ve m_1\rVert^2 \upto \lVert\ve s_4 - \ve
  m_4\rVert^2\bigr)$ depends polynomially on the
  coefficients~$a_{i,j}$, $b_i$, and~$c_i$ of $A$, $S$, and $\Mini$,
  respectively, so there is a polynomial
  $F(x_{1,1} \upto x_{3,4},y_1 \upto y_6,z_1 \upto z_5)$ in~$23$
  indeterminates such that
  \begin{equation} \label{2eqF}%
    f_M\bigl(\lVert\ve s_1 - \ve m_1\rVert^2 \upto \lVert\ve s_4 - \ve
    m_4\rVert^2\bigr) = F(a_{1,1} \upto a_{3,4},b_1 \upto b_6,c_1
    \upto c_5).
  \end{equation}
  Consider the matrix $X := (x_{i,j})_{1 \le i,j \le 3}$ and the ideal
  $I \subseteq \RR[x_{1,1} \upto x_{3,4}]$ (the ring of polynomials
  in~$12$ indeterminates) generated by the polynomial $\det(X) - 1$
  and by the coefficients of the matrix $X \cdot X^T - I_3$. Clearly a
  matrix $A \in \Aff(3)$ lies in $\ASO(3)$ if and only if every
  polynomial from $I$ vanishes when evaluated at the coefficients of
  $A$. This implies that $I$ is contained in the vanishing ideal of
  $\ASO(3)$, which we write as $I \subseteq \Id(\ASO(3))$. We claim
  equality. Indeed, using MAGMA~[\citenumber{magma}], we can verify
  that $I$ is a radical ideal and equidimensional of
  dimension~$6$. (In fact, MAGMA computes over $\QQ$ instead of $\RR$,
  but over fields of characteristic~$0$, the algorithms for computing
  radical ideals and equidimensional parts yield the same result when
  passing to a field extension, see
  \mycite[Chapter~4]{Greuel.Pfister}.) So $I \subsetneqq \Id(\ASO(3))$
  would imply that $\ASO(3)$ has dimension $< 6$. But it is well known
  that $\dim(\ASO(3)) = 6$, so indeed $I = \Id(\ASO(3))$.

  Now choose a Gr\"obner basis $G$ of $I$ with respect to an arbitrary
  monomial ordering.
  Viewing the~$x_{i,j}$ as the main indeterminates,
  we can form a normal form $\tilde{F} := \NF_G(F)$ of the polynomial
  $F$ with respect to $G$. Let
  $J \subseteq \RR[y_1 \upto y_6,z_1 \upto z_5]$ be the ideal
  generated by the coefficients of $\tilde{F}$.
  \begin{claim} \label{2Claim3}%
    If there exists $k \ge 0$ such that
    \begin{equation} \label{2eqJ}%
      (z_3 z_5)^k y_i \in J
    \end{equation}
    for $i \in \{1 \upto 6\}$, then \cref{2Claim2} and therefore
    \cref{2tFixed,2cDPWLV} follow.
  \end{claim}
  To prove the claim, let $\ve s_1 \upto \ve s_4 \in \RR^3$, not all
  equal, and let $\Mini$ be given. As shown above, we may assume that
  the~$\ve s_i$ are as in~\cref{2eqS} and $\Mini$ as
  in~\cref{2eqNewMini}. So $c_3 c_5 \ne 0$ and $b_i \ne 0$ for at
  least one~$i$. By~\cref{2eqJ}, this implies that at least one
  generator~$g$ of $J$ satisfies
  $g(b_1 \upto b_6,c_1 \upto c_5) \ne 0$. By the definition of $J$,
  this means that
  $\tilde{F}(x_{1,1} \upto x_{3,4},b_1 \upto b_6,c_1 \upto c_5) \ne
  0$. Write $F = \sum_{i=1}^m g_i t_i$ with
  $g_i \in \RR[y_1 \upto y_6,z_1 \upto z_5]$ and~$t_i$ power products
  of the~$x_{i,j}$. The linearity of the normal form map implies
  $\tilde{F} = \sum_{i=1}^m g_i \NF_G(t_i)$, so
  \begin{multline*}
    \tilde{F}(x_{1,1} \upto x_{3,4},b_1 \upto b_6,c_1 \upto c_5) =
    \sum_{i=1}^m g_i(b_1 \upto b_6,c_1 \upto c_5) \NF_G(t_i) = \\
    \NF_G\bigl(F(x_{1,1} \upto x_{3,4},b_1 \upto b_6,c_1 \upto
    c_5)\bigr).
  \end{multline*}
  Since this is non-zero and since a polynomial in
  $\RR[x_{1,1} \upto x_{3,4}]$ has normal form zero if and only if it
  lies in $I$, we obtain
  $F(x_{1,1} \upto x_{3,4},b_1 \upto b_6,c_1 \upto c_5) \notin I$.
  Since $I = \Id(\ASO(3))$, this implies that there exists
  $A \in \ASO(3)$ with coefficients $a_{1,1} \upto a_{3,4}$ such that
  \[
  F(a_{1,1} \upto a_{3,4},b_1 \upto b_6,c_1 \upto c_5) \ne 0.
  \]
  By~\cref{2eqF} this means
  $f_M\bigl(\lVert\ve s_1 - \ve m_1\rVert^2 \upto \lVert\ve s_4 - \ve
  m_4\rVert^2\bigr) \ne 0$,
  so indeed \cref{2Claim2} follows from~\cref{2eqJ}.

  It remains to show~\cref{2eqJ}, and this can be checked with the
  help of a computer. It turns out that~\cref{2eqJ} holds with
  $k = 2$. For the verification, we used MAGMA and proceeded as
  follows: 
  \begin{itemize}
  \item It is straightforward to compute the polynomials $F$ and
    $\tilde{F}$ according to their definitions. Let
    $C \subset \RR[y_1 \upto y_6,z_1 \upto z_5]$ be the set of all
    coefficients of $\tilde{F}$.
  \item Using an additional indeterminate~$t$, we computed the set
    $C^\h$ of homogenizations of the polynomials in $C$ with respect
    to~$t$.
  \item We computed a truncated Gr\"obner basis $G^\h$ of the ideal
    generated by $C^\h$ of degree~$6$.
  \item We checked that $\NF_{G^\h}\bigl(t (z_3 z_5)^2 y_i\bigr) = 0$
    for $i = 1 \upto 6$. This shows that $t (z_3 z_5)^2 y_i$ is an
    $\RR[y_1 \upto y_6,z_1 \upto z_5,t]$-linear combination of the
    polynomials in $C^\h$. Setting~$t = 1$ shows that $(z_3 z_5)^2 y_i
    \in J$.
  \end{itemize}
  The MAGMA code for running these computations, and those for the
  proof of \cref{3tMain}, is available at
  \url{https://purr.purdue.edu/publications/3105/1}.  The total
  computation time was less than a tenth of a second.
\end{proof}


\section{A drone carrying microphones and a
  loudspeaker} \label{3sCarry}%

As before, assume that we are given a room and a drone. In contrast to
the last section, assume that not only four microphones but also a
loudspeaker are mounted on the drone. Again we say that the drone is
in a \df{good position} if \cref{1aDetect} detects no walls that are
not really there. The goal of this section is to prove the following
result. It is more delicate than \cref{2tFixed} since the mirror
points~$\ve s$ are not fixed but move as the drone moves.

\begin{theorem} \label{3tMain}%
  Consider a given room, by which we understand an arrangement of
  walls, which may include ceilings, floors, and sloping walls. Also
  consider a drone that carries four non-coplanar microphones and a
  loudspeaker at fixed locations on its body. Place the drone in the
  room at a random position, which means that not only the location of
  the drone's center of gravity is chosen at random, but also its
  pitch, yaw, and roll. Then with probability~$1$ the drone is in a
  good position. More precisely, within the configuration space
  $\RR^3 \times \SO(3)$ of possible drone positions, the bad ones lie
  in a subvariety of dimension~$\le 5$.
\end{theorem}

\newcommand{\Lini}{\ve L_{\operatorname{ini}}}%
\setcounter{claim}{1}%
\renewcommand{\theclaim}{\arabic{claim}'}%
\begin{proof}
  The proof is similar to the previous one, but more complicated. We
  use some of the notation from the previous proof. In particular, the
  initial positions of the microphones are given by a matrix $\Mini$
  as~\cref{2eqMini}, but in addition we are given a vector
  $\Lini \in \RR^3$, so the actual microphone and loudspeaker
  positions are determined by
  \begin{equation} \label{3eqM}%
    M := \left(
      \begin{array}{cccc}
        \\ \ve m_1 & \ve m_2 & \ve m_3 & \ve m_4 \\ \\
        \hline 1 & 1 & 1 & 1
      \end{array}\right) = A \cdot \Mini \quad \text{and} \quad
    \left(
      \begin{array}{c}
        \\ \ve L \\ \\ \hline 1
      \end{array}\right) =
    A \cdot \left(
      \begin{array}{c}
        \\ \Lini \\ \\ \hline 1
      \end{array}\right)
  \end{equation}
  with $A \in \ASO(3)$. \Cref{2Claim1} and its proof carry over
  verbatim to the present situation. \Cref{2Claim2} has to be
  modified as follows:
  \begin{claim} \label{3Claim2}%
    Let $W_1,W_2,W_3,W_4 \subset \RR^3$ be four planes that are not
    all equal. Then there exists $A \in \ASO(3)$ such that
    \[
    f_M\bigl(\lVert\re_{W_1}(\ve L) - \ve m_1\rVert^2 \upto
    \lVert\re_{W_4}(\ve L) - \ve m_4\rVert^2\bigr) \ne 0.
    \]
    (The above expression depends on $A$ by~\cref{3eqM}.)    
  \end{claim}
  The proof that \cref{3Claim2} implies \cref{3tMain} is also as
  above. One only has to observe that in
  $f_M\bigl(\lVert\re_{W_1}(\ve L) - \ve m_1\rVert^2 \upto
  \lVert\re_{W_4}(\ve L) - \ve m_4\rVert^2\bigr)$,
  not only the~$\ve m_i$ but also $\ve L$ now depend on~$A$, but since
  the $\re_{W_i}$ are linear maps, the expression again depends
  polynomially on the coefficients of~$A$.

  For the proof of \cref{3Claim2} we must represent the walls
  $W_i$ in an appropriate way. For this, we choose normal
  vectors~$\ve w_i \in \RR^3$ such that $\ve w_i \in W_i$ if $W_i$
  does not contain the origin~$\ve 0$ of the coordinate system. So
  \begin{equation} \label{3eqWi}%
    W_i =
    \begin{cases}
      \bigl\{\ve x \in \RR^3 \mid \langle\ve w_i,\ve x\rangle =
      0\bigr\} &
      \text{if}\ \ve 0 \in W_i \\
      \bigl\{\ve x \in \RR^3 \mid \langle\ve w_i,\ve x\rangle =
      \lVert\ve w_i\rVert^2\bigr\} & \text{if}\ \ve 0 \notin W_i
    \end{cases}.
  \end{equation}
  (More formally, $W_i$ is given by~$\ve w_i$ and ``true'' or
  ``false'' indicating which of the above formulas is to be used.) The
  reflection of~$\ve L$ at $W_i$ is
  \begin{equation} \label{3eqRef}%
    \re_{W_i}(\ve L) =
    \begin{cases}
      \ve L - 2 \alpha_i \ve w_i & \text{if} \ \ve 0 \in W_i \\
      \ve L + 2 (1 - \alpha_i) \ve w_i & \text{if} \ \ve 0 \notin W_i
    \end{cases}
    \quad \text{with} \quad \alpha_i := \frac{\langle\ve w_i,\ve
      L\rangle}{\lVert\ve w_i\rVert^2}.
  \end{equation}
  Form the matrix
  $W = \begin{pmatrix} \ve w_1 & \ve w_2 & \ve w_3 & \ve
    w_4 \end{pmatrix} \in \RR^{3 \times 4}$ and set~$r :=
  \rank(W)$. We reorder the~$\ve w_i$ such
  that~$\ve w_1 \upto \ve w_r$ become linearly independent. Since we
  wish to prove \cref{3Claim2}, we also need to reorder
  the~$\ve m_i$ and therefore the columns of the matrices $M$ and
  $\Mini$. From the definition~\cref{1eqF} of~$f_M$ we see that this
  does not change
  $f_M\bigl(\lVert\re_{W_1}(\ve L) - \ve m_1\rVert^2 \upto
  \lVert\re_{W_4}(\ve L) - \ve m_4\rVert^2\bigr)$.  Computing the
  intersection of the first~$r$ walls amounts to solving a linear
  system of rank~$r$ with~$r$ equations, given by~\cref{3eqWi}. Since
  this is solvable, we may choose the origin~$\ve 0$ of the coordinate
  system such that $\ve 0 \in W_i$ for $i \le r$. So the number~$l$ of
  walls $W_i$ with $\ve 0 \in W_i$ satisfies $r \le l \le 4$, and we
  may reorder the walls again such that $\ve 0 \in W_i$ if and only if
  $i \le l$. By the hypothesis in \cref{3Claim2} that not all
  walls are equal, the case $r = 1$ and $l = 4$ can be excluded. As in
  the previous proof we use QR-decomposition and thus assume $W$ to be
  upper triangular. Its first~$r$ diagonal entries are therefore
  non-zero. Since $\ve 0 \in W_i$ for $i \le r$ we may rescale the
  fist~$r$ normal vectors~$\ve w_i$ such that the these diagonal
  entries become~$1$. Moreover, since $\rank(W) = r$, only the
  first~$r$ rows of $W$ can be non-zero. In summary, we obtain
  \begin{equation} \label{3eqW}%
    \begin{split}
      W =
      \begin{pmatrix}
        1 & b_1 & b_2 & b_3 \\
        0 & 1 & b_3 & b_5 \\
        0 & 0 & 1 & b_6
      \end{pmatrix} \ \text{if} \ r = 3,\quad & W =
      \begin{pmatrix}
        1 & b_1 & b_2 & b_3 \\
        0 & 1 & b_4 & b_5 \\
        0 & 0 & 0 & 0
      \end{pmatrix} \ \text{if} \ r = 2, \quad \text{and} \\
      W & =
      \begin{pmatrix}
        1 & b_1 & b_2 & b_3 \\
        0 & 0 & 0 & 0 \\
        0 & 0 & 0 & 0
      \end{pmatrix} \ \text{if} \ r = 1, \quad \text{where always} \
      b_i \in \RR.
    \end{split}
  \end{equation}
  Having reordered the $W_i$ and the columns of $M$ and $\Mini$, and
  having chosen a suitable Cartesian coordinate system, we can modify
  $\Mini$ as in the previous proof. (This includes a rescaling of the
  coordinate system, after which the above rescaling of the $\ve w_i$
  can be done.) So $\Mini$ is given by~\cref{2eqNewMini}. Moreover,
  we write $\Lini = (c_6,c_7,c_8)^T$.

  Set $\eta := \prod_{i=1}^4 \lVert\ve w_i\rVert^2$, which is non-zero
  and serves as a common denominator for the~$\alpha_i$
  in~\cref{3eqRef} and also for the
  $\lVert\re_{W_i}(\ve L) - \ve m_i\rVert^2$. Therefore
  $\eta^2 f_M\bigl(\lVert\re_{W_1}(\ve L) - \ve m_1\rVert^2 \upto
  \lVert\re_{W_4}(\ve L) - \ve m_4\rVert^2\bigr)$
  depends polynomially on the coefficients~$a_{i,j}$, $b_i$, and~$c_i$
  of $A$, $W$, and $\Mini$, respectively, so there are polynomials
  $F(x_{1,1} \upto x_{3,4},y_1 \upto y_6,z_1 \upto z_8)$ in~$26$
  indeterminates and $N(y_1 \upto y_6)$ in~$6$ indeterminates such
  that
  \begin{equation} \label{3eqN}%
    \eta = N(b_1 \upto b_6)
  \end{equation}
  and
  \begin{equation} \label{3eqF}%
    \eta^2 f_M \bigl(\lVert\re_{W_1}(\ve L) - \ve m_1\rVert^2 \upto
    \lVert\re_{W_4}(\ve L) - \ve m_4\rVert^2\bigr) = F(a_{1,1} \upto
    a_{3,4},b_1 \upto b_6,c_1 \upto c_8)
  \end{equation}
  In reality,
  $f_M\bigl(\lVert\re_{W_1}(\ve L) - \ve m_1\rVert^2 \upto
  \lVert\re_{W_4}(\ve L) - \ve m_4\rVert^2\bigr)$ and~$\eta$ also
  depend on~$r$ and~$l$ since the~$\ve w_i$ and $\re_{W_i}(\ve L)$ do,
  according to~\cref{3eqW} and~\cref{3eqRef}. To express these
  dependencies, we write $F_{r,l}$ and $N_r$ instead of $F$ and
  $N$. As before, we consider the ideal
  $I = \Id\bigl(\ASO(3)\bigr) \subset \RR[x_{1,1} \upto x_{3,4}]$ and
  the normal form $\tilde{F}_{r,l} = \NF_G(F_{r,l})$ with respect to a
  Gr\"obner basis $G$ of $I$. We also write
  $J_{r,l} \subseteq \RR[y_1 \upto y_6,z_1 \upto z_8]$ for the ideal
  generated by the coefficients of $\tilde{F}_{r,l}$.
  \begin{claim} \label{3Claim3}%
    If for every $1 \le r \le 3$ and $r \le l \le 4$ with $l - r < 3$
    there exist non-negative integers~$k_1,k_2$ such that
    \begin{equation} \label{3eqJ}%
      (z_3 z_5)^{k_1} N^{k_2} \in J_{r,l},
    \end{equation}
    then \cref{3Claim2} and therefore \cref{3tMain} follow.
  \end{claim}
  The proof is almost identical to the one of \cref{2Claim3}. Let
  $W_i$ be planes as in \cref{3Claim2}, and let $\Mini$ also be
  given. After reordering the $W_i$ and the columns of $\Mini$ and
  after choosing a convenient coordinate system, we obtain~$r$ and~$l$
  as in \cref{3Claim3} such that the arguments of
  $f_M\bigl(\lVert\re_{W_1}(\ve L) - \ve m_1\rVert^2 \upto
  \lVert\re_{W_4}(\ve L) - \ve m_4\rVert^2\bigr)$ are given
  by~\cref{3eqRef}, \cref{3eqW}, \cref{3eqM},
  and~\cref{2eqNewMini}. We have $c_3 c_5 \eta \ne 0$,
  so~\cref{3eqJ} together with~\cref{3eqN} implies
  $g(b_1 \upto b_6,c_1 \upto c_8) \ne 0$ for some generator~$g$ of
  $J$. Precisely as in the proof of \cref{2Claim3} this shows the
  existence of $A \in \ASO(3)$ with coefficients~$a_{i,j}$ such that
  $F(a_{1,1} \upto a_{3,4},b_1 \upto b_6,c_1 \upto c_8) \ne 0$.
  By~\cref{3eqF} this implies
  $f_M\bigl(\lVert\re_{W_1}(\ve L) - \ve m_1\rVert^2 \upto
  \lVert\re_{W_4}(\ve L) - \ve m_4\rVert^2\bigr) \ne 0$, which is the
  assertion of \cref{3Claim2}.

  What is left is to show~\cref{3eqJ}, for which we proceeded as in
  the proof of~\cref{2eqJ}. Here it turns out that for each~$r$
  and~$l$, a truncated Gr\"obner basis $G^\h$ (with the notation of
  the previous proof) of degree~$16$ suffices to show that
  $\NF_{G^\h}\bigl((z_3 z_5 N)^2\bigr) = 0$. The total computation
  time was roughly five minutes on a workstation, highlighting that
  the case of a loudspeaker carried by a drone is much more difficult
  than the case of a loudspeaker at a fixed position.
\end{proof}

Note that although it took roughly five minutes of computation time to
verify \cref{3tMain}, the actual wall detection algorithm performed by
the drone only requires evaluating simple expressions such as the one
given in \cref{1exStandard}, which is extremely fast.

%
%

\section{Conclusion}
We have shown that the problem of reconstructing an arrangement of
walls from the first order echoes of a single sound impulse acquired
by 4 microphones on a drone is generically well-posed.  The first
order echoes provide us with a list of distances from the microphones
to the walls. We assume that we do not know which distance goes with
which wall.  Both the case where the speaker producing the sound is
fixed in the room, and the case where it is carried on the drone, were
considered.  Specifically, our results show that, when the drone is in
a generic position and orientation, one can obtain four points on each
wall that are heard by all four microphones, and these walls are
guaranteed to exist (no ghost walls).  Our set up assumes exact
measurements and infinite precision calculations.  In future work, we
plan to study the more practical problem of reconstructing the walls
from noisy distance measurements.

While the formulation of our problem focuses on microphones mounted on
a drone, our results apply to many other application scenarios. For
example, the microphones and loudspeaker could be mounted on a car
moving on the road, a robot navigating in an indoor environment, or an
underwater vehicle exploring a wreck in the ocean.  Some of these
situations put restrictions on rotations and translations that can be
applied to the microphone configuration.  The impact of such
restrictions on the reconstruction problem will also be studied in
future work.

\end{document}